\documentclass[12pt]{article}

\usepackage[T1]{fontenc}
\usepackage[latin1]{inputenc}
\usepackage{bbm}
\usepackage{stmaryrd}
\usepackage{latexsym}
\usepackage{amsfonts}
\usepackage{amsmath,amssymb,amscd}
\usepackage{yfonts}
\usepackage{dsfont}
\usepackage{mathabx}
\usepackage[dvips]{graphicx}
\usepackage{epsfig}
\usepackage{psfrag}
\usepackage[font={small}]{caption}
\usepackage{color}
\usepackage{float}
\usepackage{upgreek}
\usepackage{tikz}
\usepackage{tikz-cd}
\usepackage{relsize}
\usepackage{url}

\DeclareFontFamily{U}{txsyc}{}
\DeclareFontShape{U}{txsyc}{m}{n}{
   <-> txsyc%
}{}
\DeclareFontShape{U}{txsyc}{bx}{n}{
   <-> txbsyc%
}{}
\DeclareFontShape{U}{txsyc}{l}{n}{<->ssub * txsyc/m/n}{}
\DeclareFontShape{U}{txsyc}{b}{n}{<->ssub * txsyc/bx/n}{}
\DeclareSymbolFont{symbolsC}{U}{txsyc}{m}{n}
\SetSymbolFont{symbolsC}{bold}{U}{txsyc}{bx}{n}
\DeclareFontSubstitution{U}{txsyc}{m}{n}
\DeclareMathSymbol{\df}{\mathrel}{symbolsC}{"42}
\DeclareMathSymbol{\fd}{\mathrel}{symbolsC}{"43}
\DeclareMathSymbol{\lJoin}{\mathrel}{symbolsC}{"58}
\DeclareMathSymbol{\rJoin}{\mathrel}{symbolsC}{"59}

\newcommand{\f}[2]{\frac{#1}{#2}}
\newcommand{\cA}{\mathcal{A}}

\newcommand{\cF}{\mathcal{F}}

\newcommand{\cL}{\mathcal{L}}

\newcommand{\cP}{\mathcal{P}}

\newcommand{\cS}{\mathcal{S}}

\newcommand{\EE}{\mathbb{E}}

\newcommand{\LL}{\mathbb{L}}
\newcommand{\NN}{\mathbb{N}}
\newcommand{\PP}{\mathbb{P}}

\newcommand{\RR}{\mathbb{R}}

\newcommand{\ZZ}{\mathbb{Z}}

\newcommand{\fs}{\mathfrak{s}}

\newcommand{\iy}{\infty}

\newcommand{\lt}{\left}
\newcommand{\me}{\medskip}

\newcommand{\ri}{\rightarrow}
\newcommand{\rt}{\right}
\newcommand{\sm}{\smallskip}

\newcommand{\wi}{\widetilde}
\newcommand{\wit}{\widehat}

\newcommand{\fo}{\forall\ }

\newcommand{\lve}{\lt\vert}
\newcommand{\lVe}{\lt\Vert}

\newcommand{\rve}{\rt\vert}
\newcommand{\rVe}{\rt\Vert}

\newcommand{\st}{\,:\,}

\newcommand{\un}{\mathds{1}}

\newcommand{\bq}{\begin{eqnarray*}}
\newcommand{\bqn}[1]{\begin{eqnarray}\label{#1}}
\newcommand{\eq}{\end{eqnarray*}}
\newcommand{\eqn}{\end{eqnarray}}
\newcommand{\wwtbp}{\hfill $\blacksquare$\par\me\noindent}

\newcommand{\thistitlepagestyle}{}
\newcommand{\lin}{\llbracket}
\newcommand{\rin}{\rrbracket}

\newcommand{\ttsim}{\raise.17ex\hbox{$\scriptstyle\mathtt{\sim}$}}

\newcommand{\kh}{\kern .08em}

\newtheorem{pro}{Proposition} 
\newtheorem{cor}[pro]{Corollary}
\newtheorem{lem}[pro]{Lemma}
\newtheorem{theo}[pro]{Theorem}

\renewcommand{\thepro}{\arabic{pro}}

\newenvironment{rem}
{\par\me\refstepcounter{pro}\noindent{\bf Remark \thepro\ }}
{\hfill $\square$\par\sm\noindent}

\newcommand{\proof}{\par\me\noindent\textbf{Proof}\par\sm\noindent}

\newcommand{\tv}{_\mathrm{tv}}

\title{On a Markov construction of couplings}
 \author{Persi Diaconis${}^\dagger$ and Laurent Miclo${}^\ddagger$}
 \date{\vbox{\copy0
 \vskip5mm
 \copy1
}
 }

\begin{document}

\setbox0=\vbox{
\large
\begin{center}
${}^\dagger$ Department of Mathematics\\ 
Department of Statistics\\
Stanford University
\end{center}
}
\setbox1=\vbox{
\large
\begin{center}
${}^\ddagger$ Toulouse School of Economics\\
Institut de Mathématiques de Toulouse\\
CNRS and University of Toulouse
\end{center}
}
\setbox6=\vbox{
\hbox{miclo@math.cnrs.fr\\[1mm]}
\hbox{Institut de Mathématiques de Toulouse\\}
\hbox{Université Paul Sabatier, 118, route de Narbonne\\}
\hbox{31062 Toulouse cedex 9, France\\[1mm]}
\hbox{Toulouse School of Economics,\\}
\hbox{1, Esplanade de l'université\\}
\hbox{31080 Toulouse cedex 06, France\\}
}
 \setbox5=\vbox{
\hbox{diaconis@math.stanford.edu\\[1mm]}
 \hbox{Department of Mathematics\\}
 \hbox{Department of Statistics\\}
 \hbox{Stanford University, USA} 
 }

\maketitle
\thistitlepagestyle
\abstract{For $N\in\NN$, let $\pi_N$ be the law of the number of fixed points of a  random permutation of $\{1, 2, ..., N\}$. Let $\cP$ be a Poisson law of parameter 1.
A classical result shows that $\pi_N$ converges to $\cP$ for large $N$ and indeed in total variation
\bq
\lVe \pi_N-\cP\rVe_{\mathrm{tv}}&\leq & \f{2^N}{(N+1)!}\eq
This implies that $\pi_N$ and $\cP$ can be coupled to at least this accuracy. This paper constructs such a coupling (a long open problem) using the machinery of intertwining of two Markov chains. This method shows promise for related problems of random matrix theory.
}
\vfill\null
{\small
}\par

\newpage 

\section{Introduction}

For $N\in\NN\df\{1,2, ...\}$, let $\pi_N$ be the law of the number of fixed points of a  random permutation of $\{1, 2, ..., N\}$. Let $\cP$ be the Poisson law of parameter 1.
A classical result, see de Montmort \cite{zbMATH03748105}, shows that $\pi_N$ converges to $\cP$ for large $N$. Indeed it is well-known (and estimates of the same order are proved  below) that in total variation
\bqn{theo1}
\f{N}{N+2} \f{2^{N+1}}{(N+1)!}\ \leq \ \lVe \pi_N-\cP\rVe\tv\ \leq \ \f{2^{N+1}-1}{(N+1)!}\eqn
\par
The total variation distance can be realised by a coupling of $\pi_N$ and $\cP$, see e.g.\ Proposition 4.7 of Levin,  Peres and Wilmer \cite{MR2466937},
and it has been a long open problem to give an explicit realization of such a coupling.
The super-exponential errors bounds in \eqref{theo1} occur in other problems such as the number of $k$-cycles in a random permutation, which has a limiting Poisson distribution of parameter $1/k$ with super-exponential error. Similar results hold for the trace of powers of random matrices for the compact classical groups $O_N$, $U_N$ and $SP_{2N}$, see e.g.\ Courteaut, Johansson and Lambert \cite{courteaut2022berryesseen}. The method introduced here shows promise for finding couplings for these problems.
For a history of Montmort's theorem, see Takacs \cite{zbMATH03675825}. For extensions and a recent literature review, see Diaconis and Fulman and Guralnick
\cite{zbMATH05566029}. At the end of this introduction we will present several attempts, successful as well as unsuccessful, to get a proof by coupling of \eqref{theo1}.
\par\me
To present our approach, for any $N\in\NN$ and any permutation $\sigma$ in the symmetric group $\cS_N$, denote $\eta_1(\sigma)$ the number of fixed point of $\sigma$:
\bq
\eta_1(\sigma)&\df& \vert \{x\in\lin N\rin\st \sigma(x)=x\}\vert
\eq
(where $\lin N\rin\df\{1, 2, ..., N\}$ and more generally, for any $n\leq n'\in\ZZ_+\df\{0,1,2, ...\}$, we write $\lin n,n'\rin\df \{n, n+1, ..., n'\}$). The number $\eta_2(\sigma)$ of 2-cyles of $\sigma$ will also play an important role:
\bq
\eta_2(\sigma)&\df&\vert\{(x,y)\in \lin N\rin^2\st x<y,\, \sigma(x)=y\hbox{ and } \sigma(y)=x\}\vert\eq
\par
Let $\nu_N$ stands for the uniform distribution on $\cS_N$, so that $\pi_N$ is its image by $\eta_1$ on $\ZZ_+$. To simplify the notation, we will often drop the exponent $N$ when referring to these probability measures.
As mentioned in \eqref{theo1}, the fixed-point law
 $\pi$ is very close to the Poisson distribution $\cP$. The bounds in \eqref{theo1} are for instance recorded in
  (1.11) page 15 of Arratia, Barbour and Tavaré \cite{zbMATH02042126} and are deduced from computations of David and Barton \cite{MR0155371}
using properties of alternating series with decreasing terms coming from the following traditional facts.\par
We have
\bqn{piD}
\fo x\in\lin 0, N\rin,\qquad \pi(x)&=&\f{D_{N-x}}{(N-x)!}\f1{x!}\eqn
where for any $n\in\ZZ_+$, $D_n$ stands for the number of derangements from $\cS_n$, namely the permutations of $\cS_n$ without fixed point (with the convention that $D_0=1$).  The formula due to 
de Montmort \cite{zbMATH03748105} gives the number of derangements:
\bqn{der}
\fo n\in\NN,\qquad D_n&=&n!\sum_{k=0}^n\f{(-1)^k}{k!}\eqn
leading to the explicit formula:
\bqn{pi}
\fo x\in\lin 0, N\rin,\qquad \pi(x)&=&\f1{x!}\sum_{k=0}^{N-x}\f{(-1)^k}{k!}\eqn\par\me
As announced above, our purpose is to  deduce  bounds on $ \lVe \pi-\cP\rVe\tv$, of the same logarithmic order as that of  \eqref{theo1}.  
Here is a sketch of the proof. We use a random transposition to construct a Markov chain on the symmetric group $\cS_N$. Then the intertwining-lumping procedure presented in Section \ref{general} and some fiddling around is used to construct a monotone birth-and-death chain with the fixed point distribution $\pi$ as reversible distribution. A similar construction gives a monotone birth-and-death chain with a Poisson stationary distribution. Revisiting next the classical coupling of two monotone birth-and-death chains leads to our bound. In more detail the intertwining-lumping construction produces  the 
penta-diagonal Markov kernel $P$ on
\bq
V&\df&\lin 0, N-2\rin\sqcup\{N\}\eq
given by
\bqn{P}
\fo x\in V,\qquad
\lt\{\begin{array}{rcl}
P(x,x-1)&=&\f{x(N-x)}{N(N-1)}\\
P(x,x-2)&=&\f{x(x-1)}{N(N-1)}\\
P(x,x+1)&=&\f{N-x-2p(x)}{N(N-1)}\\
P(x,x+2)&=&\f{2p(x)}{N(N-1)}\\
P(x,x)&=&1-P(x,x-1)-P(x,x-2)-P(x,x+1)-P(x,x+2)
\end{array}\rt.
\eqn
where \bqn{p}
x\in V,\qquad p(x)&\df &\EE_\nu[\eta_2\vert \eta_1=x]\eqn
(the conditional expectation is with respect to the uniform measure $\nu$ on $\cS_N$). As explained in Section~\ref{general} below, this chain is a projection of Markov chains on conjugacy classes derived from multiplication from random transpositions.\par
Note that $P$ does not allow to get out of $V$: we have $P(0,-1)=P(0,-2)=P(1,-1)=P(N, N-1)=0$ and
$P(N-3,N-1)=P(N-2,N-1)=P(N,N+1)=P(N,N+2)=0$. For the latter equalities, we need the following observations about $p$:
obviously we have $p(N)=0$ and 
the value $p(N-2)$ is 1, since knowing that $\eta_1=N-2$, we necessarily have $\eta_2=1$.
Similarly, the value $p(N-3)$ is 0, since knowing that $\eta_1=N-3$, we necessarily have $\eta_2=0$ (and the number of 3-cycles is equal to 1).
\par
\sm
By our construction,  the probability measure $\pi$ will naturally appear to be reversible for the Markov kernel $P$. Furthermore the reversibility of $P$ (without even knowing the reversible probability) in conjunction with $p(N)=p(N-3)=0$
and $p(N-2)=1$ are sufficient to determine $P$ and by consequence the other values of $p$ and those of $\pi$.
These features can be translated into convenient estimates on $p$, leading to quantitative couplings of the Markov chains whose transitions are dictated by $P$ with other Markov chains whose invariant measure is the conditioning of $\cP$ on $V$ (more conveniently, we will restrict our attention to the state space $\lin 0, N-4\rin\subset V$). These bounds are carried out in Section \eqref{priori} and we will deduce the convergence
\bqn{lim}
\lim_{N\ri\iy} \f1{N\ln(N)}\ln(\lVe \pi-\cP\rVe_{\mathrm{tv}})&= &-1\eqn
of the right logarithmic order.
\par\me
Let us now list several attempts to prove \eqref{theo1} via coupling arguments, as well as some remarks.

\subsection{A failed effort}

This section records a natural coupling, indeed one that extends to all the classical compact groups and their Weyl groups. By the law ``natural yields right'', this should work to give good error bounds, alas it doesn't! \par
Let $(X_n)_{n\in\NN}$ be independent $\{0,1\}$-valued random variables with \bq
\fo n\in\NN,\qquad \PP[X_n=1]&=&\f1n\ =\ 1-\PP[X_n=0]\eq\par
Define for all $N\in\NN$,
\bq
S_N&\df& X_1X_2+X_2X_3+\cdots+X_{n-1}X_N+X_N\\
S_\iy&\df& X_1X_2+X_2X_3+\cdots
\eq
In the unpublished paper of Diaconis and Mallows \cite{Diaconis_new}, recorded in
Diaconis and Forrester \cite{zbMATH06699672}, it is shown that for any $ k\in \ZZ_+$,
\bq
 \PP[S_N=k]\ =\ \pi_N(k)&\hbox{and}& \PP[S_\iy=k]\ =\ \cP(k)
 \eq
\par Thus the joint law of $(S_N,S_\iy)$ makes a natural coupling. Alas, $S_\iy-S_N=X_N(X_{N+1}-1)+X_{N+1}X_{N+2}+X_{N+2}X_{N+3}+\cdots $ has typical distance of order $1/N$. For more background and details, see Diaconis and Forrester \cite{zbMATH06699672}.

\subsection{A successful and strange coupling from computer science}

Jim Pitman has explained a fascinating construction of a super exponential coupling due to computer scientists
Duchon and Duvignau \cite{zbMATH06928911} and Duchon and Duvignau \cite{zbMATH06661374}. Pitman's development of these ideas is unpublished \cite{Pitman_triangle}. We thank him for permission to state his results.
The construction calls for a countable collection $(U_n)_{n\in\NN}$ of independent random variables uniformly distributed on $[0,1]$. Define
\bq
S&\df& \min\{n\geq 1\st U_n<U_{n+1}\},\qquad \hbox{time of first ascent}\\
T&\df& \min\{n\geq 2\st U_n>\max(U_{n-1},U_{n+1})\},\qquad \hbox{time of first peak}\\
M&\df& S-\delta_{T-S\mathrm{\, is\, odd}}\eq
\begin{theo} The distribution of the random variable $M$ is the Poisson law of parameter 1.
\end{theo}
Define further for fixed $N\in\NN$,
\bq
S_N&\df& \min(S,N)\\
T&\df& \min(T,N)\\
M_N&\df& S_N-\delta_{T_N-S_N\mathrm{\, is\, odd}}\eq
\begin{theo} The random variable $M_N$ has the law of the number of fixed points of a  random permutation of $\lin N\rin$. 
\end{theo}
As a consequence of the two previous theorems, we get
\begin{cor}
For any $N\in\NN$, we have
\bq
\lVe \cL(M)-\cL(M_n)\rVe_{\mathrm{tv}}&\leq & \PP[T>N]\ \leq \ \f{2^N}{(N+1)!}\eq
\end{cor}
\par
This result seems magical and the present paper records an effort to find a proof using more standard tools which might permit generalization. We also hope to study it on its own at least to generalize to the law of the number of $k$-cycles.

\subsection{Unstability of the super-exponential bounds}

The previous super-exponential bounds are delicate. Consider for example the number of fixed points in the first $N-1$ places of a   random permutation of $\lin N\rin$. 
This quantity too has an approximate Poisson distribution of parameter 1 but the total variation distance between these two laws is of order $1/N$.
\par
Similarly, for any $\theta\in (0,1)$, the number of fixed points in places $ \lin\lfloor \theta N\rfloor\rin$ has a Poisson law of parameter $\theta$ as limiting law. Indeed the point process on $[0,1]$ which has an event at $k/N$ if and only if a random permutation $\sigma$ of $\cS_N$ satisfies $\sigma(k)=k$ is well approximated by a unit rate Poisson process. But these approximations are only accurate up to order $1/N$.

\subsection{Equality of first $N$ moments}

For $N\in\NN$, consider two random variables $X_N$ and $X_\iy$ respectively distributed according to $\pi_N$ and $\cP$. The high order of contact between these two laws can be captured by moments. Indeed Diaconis and Shahshahani \cite{MR1274717} show
\bq
\fo k\in\lin 0, N\rin,\qquad \EE[X_N^k]&=&\EE[X_\iy^k]\eq
\par
Similar results hold for the joint mixed moments of the number of $k$-cycles and for compact classical groups.

\subsection{The Markov approach}

It is related to Stein's method, see e.g.\ Diaconis and Holmes \cite{zbMATH02172364} or Section 4 of Chatterjee,  Diaconis and Meckes \cite{zbMATH05728575}, but the underlying philosophy is quite old.
Assume we would like to investigate some features of a given probability measure $\pi$. The Markov approach consists in introducing and studying a Markov process (in continuous time) or chain (in discrete time) encapsulating  the ``relevant characteristics'' of the underlying state space and admitting $\pi$ as invariant probability (sometimes no effort is required in this introduction, as $\pi$ is already defined as an invariant probability). 
An example of this situation is the investigation of absence of phase transition, exponential decay of correlations, or analyticity of correlation of Gibbs measures, which was done via the use of stochastic Ising processes leaving these Gibbs measures invariant, see Holley and Stroock \cite{MR428984} or Chapter 4 of the book of Liggett \cite{MR2108619}.
Our goal here is to give a new illustration of  this Markov approach by recovering the right order of \eqref{theo1}.
\par\me
The plan of the paper is as follows. In the next section we present a general procedure producing a Markov chain by projection of another Markov chain.
Reversibility is preserved by such projections. In Section \ref{penta},  the  transposition random walk on $\cS_N$ is projected in this way through $\eta_1$ to get the Markov kernel $P$ on $V$. In Section \ref{priori}, we deduce the a priori bounds on $p$ that are applied in Section~\ref{coupling} to control our couplings of Markov chains, leading to desired upper bound on the approximation of $\pi$ by $\cP$.
 In a spirit similar to that of Section \ref{priori}, in Appendix \ref{Api}, we directly recover \eqref{pi}, giving an alternative proof to the  classical inclusion-exclusion argument.
 In Appendix \ref{complements}, some complements are given about the conditional expectation $p$.

 \bigskip
 \par\hskip5mm\textbf{\large Acknowledgments:}
We thank Jim Pitman for telling us about Section 1.2. Diaconis is funded by  NSF grant 1954042. Miclo is funded by grants ANR-17-EURE-0010 and AFOSR-22IOE016.
 
 \par\sm\noindent 

\section{Projections of Markov chains}\label{general}

We present in this short section a general procedure of projection of Markov chains. We will restrict our attention  to finite state spaces to simplify the exposition and since latter we will work only with such sets, but the underlying principle is much more general.\par\me

Consider a Markov chain $(X,Y)\df(X_n,Y_n)_{n\in\ZZ_+}$ taking values in a product state space $V\times W$.
Assume that $V$ and $W$ are finite and that the transition matrix $Q$ of $(X,Y)$ is irreducible. Denote by $\mu$ its invariant measure.
Consider $r_1\st V\times W\ri V$ and $r_2\st V\times W\ri W$ the canonical projections and let $\mu_1\df r_1(\mu)$ 
 be the first 
 marginal distribution of $\mu$. Denote by $\mu_{1,2}$ the Markov kernel from $V$ to $W$ corresponding to the conditional distribution  of $r_2$ knowing $r_1$.
So we have the decomposition\
\bq
\fo (x,y)\in V\times W,\qquad \mu(x,y)&=&\mu_1(x)\mu_{1,2}(x,y)\eq
\par
Consider the Markov kernel $P$ given on $V$ via
\bq
\fo x,x'\in V\qquad P(x,x')&\df& \sum_{y,y'\in W} Q((x,y),(x',y'))\, \mu_{1,2}(x,y)\eq
and the Markov kernel $\Lambda$ from $V\times W$ to $V$ given by
\bq
\fo (x,y)\in V\times W,\,\fo x'\in V,\qquad \Lambda((x,y),x')&\df& \mu_1(x')\eq
\par
\begin{lem}\label{intert}
We have the intertwining relation
\bq
Q\Lambda&=&\Lambda P\eq
\end{lem}
\proof
On one hand, $\Lambda$ can be identified with $\mu_1$, so that $Q\Lambda=\mu_1$.\par
On the other hand, we have for any $(x,y)\in V\times W$ and $x'\in V$,
\bq
\Lambda P((x,y),x')&=&\mu_1P(x')\\
&=&\sum_{x''\in V}\mu_1(x'')P(x'',x')\\
&=&\sum_{x''\in V}\mu_1(x'') \sum_{y',y''\in W} Q((x'',y''),(x',y'))\, \mu_{1,2}(x'',y'')\\
&=&\sum_{y'\in W}\sum_{(x'',y'')\in V\times W} \mu(x'',y'')Q((x'',y''),(x',y'))\\
&=&\sum_{y'\in W}\mu(x',y')\\
&=&\mu_1(x')\eq
namely
\bq
\Lambda P&=&\mu_1\ =\ Q\Lambda\eq
\wwtbp
\par
In particular, $\mu\Lambda$ is invariant is for $P$, i.e.\ $\mu_1$ is invariant for $P$ (in fact this is just the above proof).
We also have:
\begin{lem}\label{rev}
Assume that $\mu$ is reversible for $Q$, then $\mu_1$ is reversible for $P$.
\end{lem}
\proof
Consider $f,g\in\RR^V$. We have
\bq
\mu_1[fP[g]]&=&\sum_{x\in V}\mu_1(x)f(x)\sum_{y\in W} Q[g\circ r_1](x,y)\mu_{1,2}(x,y)\\
&=&\sum_{(x,y)\in V\times W}f\circ r_1(x,y) Q[g\circ r_1](x,y) \mu(x,y)\\
&=&\mu[f\circ r_1 Q[g\circ r_1]]\\
&=&\mu[g\circ r_1 Q[f\circ r_1]]\\
&=&\mu_1[gP[f]]\eq
\wwtbp
\par
The construction above corresponds to a lumping procedure. More generally, let $W$ be a finite (or denumerable) set and $(Q(w,w'))_{w,w'\in W}$ be a Markov kernel on $W$
admitting $(\pi(w))_{w\in W}$ as stationary distribution. Given a partition of the state space $W=\bigsqcup_{v\in V} A_v$ into non-empty subsets,
reporting which $A_v$ contains the current state of a Markov chain associated to $Q$ gives a ``lumped process''. As is well-known, see e.g.\ Theorem 6.3.2 of Kemeny and Snell \cite{zbMATH03854137} or Pang \cite{zbMATH07120199}, this may not be a Markov chain. The analogous projected Markov kernel $(P(v,v'))_{v,v'\in V}$ on $V$ can defined as
\bq
\fo v,v'\in V,\qquad P(v,v')&\df& \sum_{w\in A_v,\,w'\in A_{v'}}\f{\pi(w)}{\pi(A_v)}Q(w,w')\eq
\par
Arguing as above, the probability measure $(\pi(A_v))_{v\in V}$ is invariant for $P$ (and reversible when $\pi$ is reversible for $Q$).
Defining 
\bq
\fo w\in W,\,\fo v\in V,\qquad \Lambda(w,v)&\df& \pi(A_v)
\eq
we get the intertwining relation $Q\Lambda=\Lambda P$. If the classical Dynkin condition holds, namely for any $v,v'\in V$, $Q(w,A_{v'})$ does not depend on the choice of $w\in A_v$, then the projected chain $P$ agrees with the usual lumped chain.

\section{A penta-diagonal and two birth and death Markov chains}\label{penta}

Here we apply the abstract projection scheme of the previous section in the setting of the symmetric group $\cS_N$. It is related to Chapter 12 of Stein's book \cite{zbMATH00050805}, which studies the law of the numbers of the cycles of length $l$, for all $l\in\lin N\rin$, under the uniform distribution on $\cS_N$, using a random transposition to build a reversible Markov chain.
\par\me
Consider the transposition random walk on the symmetric group $\cS_N$, whose transition matrix $T$ is given by
\bq
\fo \sigma,\,\sigma'\in\cS_N,\qquad T(\sigma,\sigma')&=&\lt\{
\begin{array}{ll}
\f2{N(N-1)}&\hbox{, if there exists a transposition $\tau$ such that $\sigma'=\tau\sigma$}\\
0&\hbox{, otherwise}
\end{array}\rt.\eq
(where permutations are seen as bijective mappings from $\lin N\rin$ and the product corresponds to the composition).
\par
The Markov kernel $T$ is reversible with respect to the uniform probability distribution $\nu$ on $\cS_N$.\par
Generalizing $\eta_1$ and $\eta_2$, for any $l\in\lin N\rin$ and $\sigma\in\cS_N$ define $\eta_l(\sigma)$ as the number of cycles of order $l$ in $\sigma$ (singleton cycles 
corresponding to fixed points).
In particular we have
\bq
\fo \sigma\in\cS_N,\qquad \eta_1(\sigma)+2\eta_2(\sigma)+\cdots +N\eta_N(\sigma)&=&N\eq
\par
Let $(\sigma(n))_{n\in\ZZ_+}$ be a Markov chain with transitions dictated by $T$ and denote
\bq
\fo n\in\ZZ_+,\qquad \eta(n)&\df& (\eta_l(\sigma(n)))_{l\in\lin N\rin}\eq
\par
It is well-known that $\eta\df(\eta(n))_{n\in\ZZ_+}$ is also a Markov chain whose transition matrix is denoted $Q$ and is reversible with respect to $\mu$ the image of $\nu$ in the mapping $\cS_N\ni\sigma\mapsto (\eta_l(\sigma))_{l\in\lin N\rin}$. Indeed this is the classical coagulation-fragmentation chain of statistical mechanics, see  Diaconis,  Mayer-Wolf,  Zeitouni and Zerner \cite{zbMATH02100738}. \par
The Markov chain $\eta$ can be written under the form $(X,Y)$ with
\bq
X&\df&\eta_1\\
Y&\df& \eta_{\lin 2,N\rin}\ \df\ (\eta_2, \eta_3,...,\eta_N)\eq
\par
We are thus in position to apply Lemmas \ref{intert} and \ref{rev}. 
\par
Our next goal is to describe the corresponding Markov kernel $P$. 
Note that the corresponding state space is
\bqn{V}
V&=&\lin 0, N\rin\setminus \{N-1\}\eqn
already met in the introduction (it is not possible for a permutation to have $N-1$ fixed points).\par\sm
Consider a permutation $\sigma\in\cS_N$. Denote $f_1, f_2, ..., f_k$ its fixed points (so that $\eta_1(\sigma)=k$) 
and let $C_1, C_2, ..., C_l$ be the other cycles of $\sigma$.
\par
Consider a transposition $\tau\fd(i,j)$. Let us describe $\eta_1(\sigma')$ with $\sigma'\df \tau\sigma $.\par
\sm
$\bullet$ If both $i$ and $j$ are fixed points of $\sigma$, then the fixed points of $\sigma'$ are the elements of $\{f_1, f_2, ..., f_k\}\setminus\{i,j\}$ 
and its non-singleton cycles are $C_1, C_2, ..., C_l$ and $(i,j)$. Thus we have $\eta_1(\sigma')=\eta_1(\sigma)-2$.\par\sm
$\bullet$ If $i$ is a fixed point of $\sigma$ and $j\in C_r$, with $r\in\lin l\rin$, then  the fixed points of $\sigma'$ are the elements of $\{f_1, f_2, ..., f_k\}\setminus\{i\}$
and its  non-singleton cycles are $C_m$ with $m\neq r$ in addition to a new cycle containing $C_r$ and $j$. Thus we have $\eta_1(\sigma')=\eta_1(\sigma)-1$.\par\sm
$\bullet$ If $i\in C_r$ and $j\in C_s$ with $r\neq s$, then  the cycles and fixed points of $\sigma'$ are the same as those of $\sigma$, except that $C_s$ and $C_r$
are merged into a new cycle. In particular we have $\eta_1(\sigma')=\eta_1(\sigma)$.\par\sm
$\bullet$ The last situation is when $i$ and $j$ belong to the same cycle $C_r$. We consider three subcases:\par\sm
- When $C_r=(i,j)$, then the fixed points of $\sigma'$ are $\{f_1, f_2, ..., f_k,i,j\}$ and its non-singleton cycles are the $C_s$, for $s\in\lin l\rin\setminus\{r\}$.
We deduce $\eta_1(\sigma')=\eta_1(\sigma)+2$.
\par\sm
- When there exists $x\in C_r$ such that $\{i,j\}=\{x,\sigma(x)\}$, assume for instance that $i=x$ and $j=\sigma(x)$ and $C_r\neq (i,j)$.
Then $i$ is a new fixed point of $\sigma'$ and its non-singleton cycles are the same as those of $\sigma$, except that the point $i$ has been removed from $C_r$. We deduce $\eta_1(\sigma')=\eta_1(\sigma)+1$.
\par\sm
- When there does not exist $x\in C_r$ such that $\{i,j\}=\{x,\sigma(x)\}$ (in particular the cardinal of $C_r$ is at least 4), then $\sigma'$ has the same fixed points as $\sigma$ and the only difference in its non-singleton cycles is that $C_r$ has been divided into two new non-singleton cycles. We deduce $\eta_1(\sigma')=\eta_1(\sigma)$.
\par\me
Integrating these observations with respect to $\tau$ uniformly distributed among all transpositions, we end up with the kernel $P$ given in \eqref{P}, with
\bq
p(x)&=&\int \eta_2\,\mu_{1,2}(x,d\eta_{\lin 2, N\rin})\eq
namely the mean of $\eta_2$ knowing $\eta_1=x$ when $\eta_{\lin 1,N\rin}$ is distributed according to $\mu$ (the above integral is in fact a sum, but the integral notation is more convenient). Note that this formulation is equivalent to \eqref{p}.
\par
The distribution $\pi$ of the number of fixed points of the uniform permutation is equal to $\mu_1$, with the notation of Section \ref{general}. According to Lemma \ref{rev}, $\pi$ is reversible for $P$, as announced in the introduction.
\par\me
In the sequel it will sometimes be more convenient to work with tri-diagonal kernels than with the penta-diagonal kernel $P$, so let us extract two 
 birth and death kernels  from $P$.
\par
The first one, denoted $\wi P$, is given by
\bqn{wiP}
\fo x\neq y\in V,\qquad
\wi P(x,y)&\df& \f1{N(N-1)}\lt\{\begin{array}{ll}
x(N-x)&\hbox{, if $x\neq N$ and $y=x-1$}\\
N-x-2p(x)&\hbox{, if $x\neq N-2$ and $y=x+1$}\\
2&\hbox{, if $x= N-2$ and $y=N$}\\
N(N-1)&\hbox{, if $x= N$ and $y=N-2$}\\
0&\hbox{, otherwise}
\end{array}\rt.\eqn
This Markov kernel is obtained by removing all transitions of the form $(x,x+2)$ and $(x,x-2)$ from $P$, except for $(N-2,N)$ and $(N,N-2)$ (because $N-1$ is not a value taken by $\eta_1$), and putting their weights to the diagonal. For the corresponding Markov chains, it amounts to forbid the jumps of size two and keep the current position instead (except for the transitions between $N-2$ and $N$).
\par
From  the fact that $\pi$ is reversible for $P$, we deduce that $\pi$ is also reversible for $\wi P$, since the property of being reversible is preserved by removing transitions (when the transitions  
in both directions along an edge are removed together). As announced,
$\wi P$ corresponds to a birth-and-death Markov transition on $V$.\par\sm
The second  birth-and-death Markov transition $\wit P$ will be useful in Appendix \ref{Api}. It is obtained by ordering
$V$ as $N-3,N-5, ...,3,  1,0,2,4, ..., N-2, N$ when $N$ is even. When $N$ is odd, rather order $V$ as $N-3,N-5, ...,2, 0,1,3, ..., N-2, N$, the following construction leads to similar results in this case, 
so let us only consider the situation where $N$ is even.
\par
Thus we define for $i\in\lin 0,N-1\rin$,
\bq
z_i&\df& \lt\{\begin{array}{ll}
2i&\hbox{, if $i\in \lin 0,N/2-2\rin$}\\[2mm]
2+2i-N&\hbox{, if $i\in \lin N/2-1, N-1\rin$}
\end{array}\rt.\eq
\par
The  Markov kernel $\wit P$ is given on $\lin 0,N-1\rin$ by
\bqn{witP}
\fo i\neq j\in V,\qquad
\wit P(i,j)&\df&\lt\{\begin{array}{ll}
P(z_i,z_j)&\hbox{, if $\vert i-j\vert=1$}\\
0&\hbox{, otherwise}
\end{array}\rt.\eqn
This construction of $\wit P$  is somewhat supplementary to that of $\wi P$: only the transitions of size two are kept, all transitions of size 1 being removed, except those between 0 and 1, to insure irreducibility.
\par
For the same reason as for $\wi P$, the kernel $\wit P$ admits $\wit \pi\df(\wit \pi(i))_{i\in\lin 0, N-1\rin}$ for reversible measure, where
\bqn{witpi}
\fo i\in\lin 0,N-1\rin,\qquad \wit \pi(i)&\df& \pi(z_i)\eqn

\section{An a priori estimate}\label{priori}

A drawback of Definition \eqref{P} of the Markov kernel $P$ is that the quantities $p(x)$, for $x\in V$, are  a priori unknown. 
We will give an explicit formula for them in Appendix \eqref{Api}, but the control of the couplings of next section only requires 
an a priori bound about them, presented in Proposition \ref{pro4.1} below.\par\me
We have seen in the introduction that $p(N-3)=p(N)=0$ and that $p(N-2)=1$.
These equalities and the fact that $P$ is reversible 
 are sufficient knowledges to deduce the following bound:
\begin{pro}\label{pro4.1}
We have for the mapping $p$ defined in \eqref{p},
\bq
\fo x\in\lin 0, N-2\rin,\qquad \vert 2p(x)-1\vert 
&\leq & \f1{(N-x-2)!}\eq
\end{pro}
\proof
Recall that Kolmogorov criterion for reversibility, see e.g.\ the book of Kelly \cite{MR2828834}, asserts that for any finite sequence $(x_0, x_1, ..., x_n)$ from $V$ with $n\in\NN$, we have
\bq
P(x_0,x_1)P(x_1,x_2)\cdots P(x_{n-1},x_n)P(x_n, x_0)&=&P(x_0,x_n)P(x_n,x_{n-1})\cdots P(x_2,x_1)\cdots P(x_{1},x_0)
\eq
\par
For given $x\in \lin 0, N-4\rin$, assuming $N\geq 4$, let us apply this formula with 
\bq
x_0&=&x\\
x_1&=&x+1\\
x_2&=&x+2\eq
\par
We get
\bq
P(x,x+1)P(x+1,x+2)P(x+2,x)&=&P(x,x+2)P(x+2,x+1)P(x+1,x)\eq
namely
\bq
\lefteqn{(N-x-2p(x))(N-x-1-2p(x+1))(x+2)(x+1)}\\&=&2p(x)(x+2)(N-x-2)(x+1)(N-x-1)\eq
i.e., since $x+1>0$,
\bq
(N-x-2p(x))(N-x-1-2p(x+1))&=&2p(x)(N-x-2)(N-x-1)
\eq
\par
To simplify notations, let us write $k(x)=2p(x)$, for any $x\in\lin 0, N-2\rin$.
The above formula is equivalent to the downward iteration, for $x\in \lin 0, N-4\rin$,
\bqn{bi}
k(x)&=&\f{(N-x)(N-x-1-k(x+1))}{(N-x-1)^2-k(x+1)}\eqn
\par
Starting from $k(N-3)=2p(N-3)=0$, we deduce iteratively $k(N-4)$, $k(N-5)$, ... down to $k(0)$.\par
For $x\in\lin 0, N-4\rin$, denote $F_x$ the rational function
\bq
\fo r\in \RR\setminus \{(N-x-1)^2\},\qquad
F_x(r)&\df& \f{(N-x)(N-x-1-r)}{(N-x-1)^2-r}\eq
so that $k(x)=F_x(k(x+1))$.
\par
For any $x\in\lin 0, N-4\rin$, 1 is a fixed point of $F_x$ (the only one in fact), since
\bq
F_x(1)&=&\f{(N-x)(N-x-1-1)}{(N-x-1)^2-1}\\
&=&1\eq
\par
Thus \eqref{bi} can be written in the convenient form
\bqn{k1}
\nonumber k(x)-1&=&F_x(k(x+1))-F_x(1)\\
&=&\int_1^{k(x+1)} F'_x(s)\,ds\eqn
which suggests  computing:
\bqn{Fp}
\fo s\in \RR\setminus \{(N-x-1)^2\},\qquad
F'_x(s)&=&-\f{(N-x)(N-x-1)(N-x-2)}{((N-x-1)^2-s)^2}
\eqn
\par
These observations lead  to a proof of the bound of Proposition \ref{pro4.1} by a backward iteration.\par
Indeed, for $x=N-2$ and $x=N-3$, the bound is true, since it is respectively implied by
\bq
2p(N-2)-1&=&2-1\\
&=&1\\
&=&\f1{0!}\\
&=&\f1{(N-2-(N-2))!}\eq
and
\bq
2p(N-3)-1&=&0-1\\
&=&-1\\
&=&-\f1{1!}\\
&=&-\f1{(N-2-(N-3))!}
\eq
\par
Consider $x\in\lin 0, N-4\rin$, we have
\bqn{NNN}
N-x\ \geq \ N-x-1\ \geq \ N-x-2\ \geq \ 2\eqn
so that $F_x'<0$. This observation and \eqref{k1} imply that if $k(x+1)>1$, then $k(x)<1$ 
and conversely, if $k(x+1)<1$, then $k(x)>1$, namely the sequence $(k(z)-1)_{z\in \lin 0, N-2\rin}$ is alternating.
\par
Let us consider separately the first case: $x=N-4$.
Since $k(N-3)=0<1$, we deduce from \eqref{Fp} with $x=N-4$, that for $s\in [k(N-3),1]=[0,1]$,
\bq
\vert F'_{N-4}(s)\vert &\leq &\f{4\times 3\times 2}{(3^2-1)^2}\ =\ \f{3}{8}\ \leq \ \f12
\eq
\par
It follows from
\eqref{k1}
that
\bq
\vert k(N-4)-1\vert&\leq &\f12\vert k(N-3)-1\vert\ \leq \ \f12\ =\ \f{1}{(N-(N-4)-2)!}\eq
\par
Let us now assume the bound of Proposition \ref{pro4.1} is true for some $x+1\in \lin 1, N-4\rin$ and let us prove it for $x$.\par
 \par
Note that we have
\bq
k(x+1)\ \leq \ 
1+\f1{(N-x-2)!}\ \leq \ 1+\f1{2!}\ \ = \ \f32\eq
so \eqref{NNN} and \eqref{Fp} imply that for $s\in [1,F_x(k(x+1))]$ (or $s\in [F_x(k(x+1)),1]$ if $F_x(k(x+1))\leq 1$),
\bq
\vert F'_x(s)\vert &\leq &\f{(N-x)(N-x-1)(N-x-2)}{((N-x-1)^2-4/3)^2}
\eq
\par
Let us show that the r.h.s.\ is bounded above by $1/(N-x-2)$. To simplify notation, write $y\df N-x-1\geq 3$, so that the desired bound amounts to
\bqn{ineg}
\f{(y+1)y(y-1)}{(y^2-4/3)^2}&\leq & \f1{y-1}\eqn
namely
\bq
(y^2-1)(y-1)y&\leq & (y^2-4/3)^2\eq
i.e.
\bq
y^4-y^3-y^2+y&\leq & y^4-\f83y^2+\f{16}9\eq
or $g(y)\geq 0$, where
\bq
\fo y\geq 3,\qquad g(y)&\df& y^3-\f53y^2-y+\f{16}9\eq\par
We compute
\bq
\fo y\geq 3,\qquad g'(y)&=&3y^2-\f{10}{3}y-1\eq
and the largest zero of the r.h.s is
\bq
\f1{18}(10+\sqrt{208})&<&3\eq
\par It follows that $g$ is increasing on $[3,+\iy)$ and we compute
\bq
g(3)&=&27-15+3+\f{16}9\ > \ 0
\eq
showing the validity of \eqref{ineg}.\par
We deduce from \eqref{k1} that 
\bq
\vert k(x)-1\vert &\leq & \lve\int_1^{k(x+1)} \f1{N-x-2}\,ds\rve\\
&\leq & \f1{N-x-2}\vert k(x+1)-1\vert\\
&\leq &  \f1{(N-x-2)!}\eq
where we took into account the iteration assumption, namely $\vert k(x+1)-1\vert\leq 1/((N-x-3)!)$.
\wwtbp
\par
\begin{rem}
The observation made after \eqref{NNN} implies more precisely  that for  $x\in \lin 0, N-2\rin$, $2p(N-2-x)-1$ is positive for even $x$ and negative for odd $x$.
\end{rem}
\par
These computations, especially the iteration relation \eqref{bi}, also show that the reversible couple $(\pi, P)$  is well-defined by $p(N-3)=0=p(N)$ and $p(N-2)=1$: no further information 
are needed for its investigation, in particular not the interpretation of $p$ as a conditional expectation on the larger space $\cS_N$. Namely we can work only on $V$.
 \par
Let us state this construction formally:
\begin{rem}\label{12}
Consider $\bar P$ defined in \eqref{P} with the $p(x)$ replaced by some $\bar p(x)\geq 0$, under the constraint that $\bar P$ is a Markov kernel on $\lin 0, N\rin$.
Add the constraints $\bar p(N-1)=\bar p(N)=0$ (in our previous case, $N-1$ becomes a transient point, with $P(N-1,N)=1/(N(N-1))$ and $P(N-1, N-1)= 1-P(N-1,N)$).
Assume furthermore that the values $\bar p(x)$ satisfy the  iteration \eqref{bi}. Thus $\bar P$ is 
 a function of $\bar p(N-2)$ and $\bar p(N-3)$.
Then taking $\bar p(N-2)=1/2$ and $\bar p(N-3)=1/2$ (implying $\bar p(x)=1/2$ for all $x\in\lin 0, N-4\rin$, due to the fact that 1 is a fixed point of $F_x$), we end up with a Markov kernel $\bar P$ which is reversible with respect to the restriction of the Poisson distribution on $\lin 0, N\rin$.
\end{rem}
\par
This observation is at the heart of the couplings presented in next section.

\section{A monotone coupling}\label{coupling}

Our purpose here is to prove by coupling an upper bound of the Poisson approximation of $\pi$, of the same logarithmic order as that of \eqref{theo1}.
It would be possible to push further the computations, but our main emphasis is placed on the method rather than on sharp estimates.\par\me
More precisely,  we want to show \eqref{lim}
by only using that $\pi$ is reversible with respect to the Markov kernel $\wi P$ defined in \eqref{wiP} and the a priori estimate given in Proposition \ref{pro4.1}.
\par
Instead of working on $V$ or $\ZZ_+$, we can restrict our attention to $\lin 0, N-4\rin$ (assuming $N\geq 4$).
Indeed, denote $\zeta$ the conditioning of $\cP$ to $\lin 0, N-4\rin$, since 
\bq
\lim_{N\ri\iy} \f1{N\ln(N)}\ln(\cP(\lin N-3,\iy\lin))&= &-1\eq
we easily deduce that
\bqn{N1}
\lim_{N\ri\iy} \f1{N\ln(N)}\ln(\lVe \zeta-\cP\rVe_{\mathrm{tv}})&= &-1\eqn
\par
Furthermore it is not difficult to see that
\bq
\lim_{N\ri\iy} \f1{N\ln(N)}\ln(\pi(\{ N-3,N-2, N\}))&= &-1\eq
since by a direct investigation, we get that
\bqn{N}
\pi(N)&=&\f1{N!}\\
\nonumber\pi(N-2)&=&\f12\f1{N!}\\
\nonumber\pi(N-3)&=&\f13\f1{N!}\eqn
\par
It follows that
\bqn{N2}\lim_{N\ri\iy} \f1{N\ln(N)}\ln(\lVe \pi-\check\pi\rVe_{\mathrm{tv}})&= &-1\eqn
where $\check \pi$ is the conditioning of $\pi$ to $\lin 0, N-4\rin$.\par
These limiting behaviors imply that \eqref{lim} amounts to
\bqn{lim2}
\limsup_{N\ri\iy} \f1{N\ln(N)}\ln(\lVe \check\pi-\zeta\rVe_{\mathrm{tv}})&\leq  &-1\eqn
\par
Indeed, on one hand, from \eqref{N} we get
\bq
\liminf_{N\ri\iy} \f1{N\ln(N)}\ln(\lVe \pi-\cP\rVe_{\mathrm{tv}})&\geq &\liminf_{N\ri\iy} \f1{N\ln(N)}\ln\lt(\f12\vert \pi(N)-\cP(N)\vert\rt)\\
&=&\liminf_{N\ri\iy} \f1{N\ln(N)}\ln\lt((1-e^{-1})/(N!)\rt)\\
&=&-1\eq
and on the other hand, from \eqref{N1}, \eqref{N2} and \eqref{lim2},
\bq
\lefteqn{\limsup_{N\ri\iy} \f1{N\ln(N)}\ln(\lVe \pi-\cP\rVe_{\mathrm{tv}})}\\&\leq  &\max \lt(\limsup_{N\ri\iy} \f1{N\ln(N)}\ln(\lVe \pi-\check\pi\rVe_{\mathrm{tv}},
\limsup_{N\ri\iy} \f1{N\ln(N)}\ln(\lVe \check\pi-\zeta\rVe_{\mathrm{tv}},\limsup_{N\ri\iy} \f1{N\ln(N)}\ln(\lVe \zeta-\cP\rVe_{\mathrm{tv}}\rt)\\
&=&-1\eq
\par
Thus it remains to prove \eqref{lim2}.\par
By reversibility of $\wi P$ with respect to $\pi$, we have that $\check \pi$ is reversible with respect to the birth and death Markov kernel  $\check P$ given by
\bq
\fo x\neq y\in \lin 0, N-4\rin,\qquad
\check P(x,y)&\df& \f1{N(N-1)}\lt\{\begin{array}{ll}
x(N-x)&\hbox{, if  $y=x-1$}\\
N-x-2p(x)&\hbox{, if $y=x+1$}\\
0&\hbox{, otherwise}
\end{array}\rt.\eq
(as usual the diagonal entries are deduced by the fact that the rows sum to 1).
\par
Consider the birth and death Markov kernel $R$ given by
\bq
\fo x\neq y\in \lin 0, N-4\rin,\qquad
R(x,y)&\df& \f1{N(N-1)}\lt\{\begin{array}{ll}
x(N-x)&\hbox{, if  $y=x-1$}\\
N-x-1&\hbox{, if $y=x+1$}\\
0&\hbox{, otherwise}
\end{array}\rt.\eq
which amounts to replacing $p(x)$ by $1/2$ in the kernel $\check P$, see Remark \ref{12} above.\par
 It is immediate to check that
$\zeta$ is reversible for $R$, since we have for any $x\in\lin 0, N-5\rin$,
\bq
\f{\zeta(x)R(x,x+1)}{\zeta(x+1)R(x+1,x)}&=&\f{(x+1)R(x,x+1)}{R(x+1,x)}\\
&=&\f{(x+1)(N-x-1)}{(x+1)(N-x-1)}\\
&=&1
\eq
\par
To simplify the notations, from now on, $\check\pi$ and $\check P$ will be written $\pi$ and $P$, we hope it will not bring confusion
with the previous $\pi$ and $P$.
\par
Consider $X\df (X(n))_{n\in\ZZ_+}$  a stationary Markov chain whose transitions are given by $P$ and whose initial law is $\pi$.
Similarly let $Y\df (Y(n))_{n\in\ZZ_+}$ be a stationary Markov chain whose transitions are given by $R$ and whose initial law is $\zeta$.
We couple them in a monotone way: namely at any time $n\in\ZZ_+$, the transition from $(X(n),Y(n))$ to $(X(n+1),Y(n+1))$ is given by sampling
an independent uniform random variable $U(n)$ on $[0,1]$ and by deciding that
\bq
X(n+1)&=&\lt\{\begin{array}{ll}
X(n)-1&\hbox{, if $U(n)<P(X(n),X(n)-1)$}\\
X(n)&\hbox{, if $P(X(n),X(n)-1)\leq U(n)<P(X(n),X(n)-1)+P(X(n),X(n))$}\\
X(n+1)&\hbox{, if $P(X(n),X(n)-1)+P(X(n),X(n))\leq U(n)$}
\end{array}\rt.\eq
and
\bq
Y(n+1)&=&\lt\{\begin{array}{ll}
Y(n)-1&\hbox{, if $U(n)<R(Y(n),Y(n)-1)$}\\
Y(n)&\hbox{, if $R(Y(n),Y(n)-1)\leq U(n)<R(Y(n),Y(n)-1)+R(Y(n),Y(n))$}\\
Y(n+1)&\hbox{, if $R(Y(n),Y(n)-1)+R(Y(n),Y(n))\leq U(n)$}
\end{array}\rt.\eq
\par
The corresponding Markov kernel on $\lin 0, N-2\rin^2$ will be denoted $S$, namely we have
\bq
\lefteqn{\hskip-40mm\fo (x,y), (x',y')\in \lin 0, N-2\rin^2, }\\
 S((x,y),(x',y'))&=&\PP[(X(n+1),Y(n+1))=(x',y')\vert (X(n),Y(n))=(x,y)]\eq\par
Consider, traditionally $\tau$ the coupling time
\bq
\tau&\df&\inf\{n\in \ZZ_+\st X(n)=Y(n)\}\eq
but also the auxiliary random chain $Z\df(Z(n))_{n\in\ZZ_+}$
\bq
\fo n\in\ZZ_+,\qquad Z(n)&\df& \sum_{k=0}^{n-1} \un_{\{X(k)=Y(k), \,X(k+1)\neq Y(k+1)\}}\eq
Their interest is that for any time $n\in\ZZ_+$, we have
\bqn{couplage2}
\nonumber\lVe \pi-\zeta\rVe_{\mathrm{tv}}&\leq & \PP[X(n)\neq Y(n)]\\
&\leq & \PP[\tau> n]+\PP [Z(n)>0]\eqn
\par
By choosing $n$ of order $N^4\ln(N)$, we will get an estimate of $\lVe \pi-\zeta\rVe_{\mathrm{tv}}$ of the order we are looking for.
\par
This resort to coupling is different from its traditional use in the quantitative investigation of convergence to equilibrium, where different lines of the same transition kernel are coupled.
The bound \eqref{couplage2} is neither good for short or long times $n$, it is interesting only for certain times, enabling us to estimate the difference between the invariant probabilities of two different transition kernels.
\par
To illustrate the difference between these approaches, let us evaluate the new term in \eqref{couplage2}:
\begin{lem}\label{lem5.4}
For any $n\in\ZZ_+$, we have
\bq
\PP [Z(n)>0]&\leq & \f{2^{N}n}{N!}\eq
\end{lem}
\proof
For any given $n\in\ZZ_+$, we have
\bqn{Zz}
\nonumber\PP [Z(n)>0]&\leq &\EE[Z(n)]\\
\nonumber&=& \sum_{k=0}^{n-1}\EE[ \un_{X(k)=Y(k), \,X(k+1)\neq Y(k+1)}]\\
&=& \sum_{k=0}^{n-1}\EE[ \un_{X(k)=Y(k)} S((X(k), Y(k)), A)]\eqn
where
\bq
A&\df& \{(x',y')\in  \lin 0, N-4\rin^2\st x'\neq y'\}
\eq
\par
Taking into account Proposition \ref{pro4.1}, we have
\bq
\fo x\in \lin 0, N-4\rin,\qquad
S((x,x),A)&\leq & \f{\vert 1-2p(x)\vert}{N(N-1)}\\
&\leq & \f1{(N-x-2)!}\f1{N(N-1)}\\
&\leq &  \f1{(N-x)!}\eq
\par
It follows that for any $k\in\lin 0, n\rin$, 
\bq
\EE[ \un_{X(k)=Y(k)} S((X(k), Y(k)), A)]&= &\sum_{x=0}^{N-4} \PP[X(k)=x=Y(k)]S((x,x),A)\\
&\leq & \sum_{x=0}^{N-4} \PP[Y(k)=x]S((x,x),A)\\
&\leq & \sum_{x=0}^{N-4}\f{1}{Z_N x!}\f1{(N-x)!}\\
&\leq &\f{1}{Z_N}\sum_{x=0}^{N}\f{1}{x!}\f1{(N-x)!}\\
&=&\f{1\times2^{N}}{Z_NN!}
\eq
where we used that $(Y(k))_{k\in\ZZ_+}$ is stationary with common distribution $\zeta$ and where 
\bq
Z_N&=&\sum_{x=0}^{N-4} \f1{x!}\ \geq \ 1\eq
\par
The desired result follows by remembering \eqref{Zz}.\wwtbp
\par
Note that the bound of the above lemma will be small even of we choose a time $n$ exponential large in $N$.
\par
In view of \eqref{couplage2} and Lemma \ref{lem5.4}, our next task is to get an estimate on $\PP[\tau>n]$ for given $n\in\ZZ_+$.
T go in this direction,
 we will need two other auxiliary random chains $\wi Z\df(\wi Z(n))_{n\in\ZZ_+}$ and $\wit Z\df(\wit Z(n))_{n\in\ZZ_+}$, defined respectively through
\bq
\fo n\in\ZZ_+,\qquad \lt\{\begin{array}{rcl}\wi Z(n)&\df& \sum_{k=0}^{n-1} \un_{X(k)\leq Y(k), \,X(k+1)> Y(k+1)}\\[2mm]
\wit Z(n)&\df& \sum_{k=0}^{n-1} \un_{X(k)\geq Y(k), \,X(k+1)< Y(k+1)}
\end{array}\rt.
\eq
as well as the hitting times of zero by $X$ and $Y$:
\bq
\tau^X_0&\df&\inf\{n\in\ZZ_+\st X(n)=0\}\\
\tau^Y_0&\df&\inf\{n\in\ZZ_+\st Y(n)=0\}
\eq
\par
Indeed, it is clear that
\bqn{tails}
\fo n\in\ZZ_+,\qquad
\PP[\tau>n]&\leq & \PP[\tau_0^X>n]+\PP[\tau_0^Y>n]+\PP[\wi Z(n)>0]+\PP[\wit Z(n)>0]\eqn
\par
It remains to estimate each of the terms of the r.h.s.\par
Let us start with the last two terms. In this respect, it is useful to remark that the Markov chain $Y$ is monotone, namely that for $x\leq y\in\lin 0, N-4\rin$, if $Y_x\df (Y_x(n))_{n\in\ZZ_+}$ and $Y_y\df (Y_y(n))_{n\in\ZZ_+}$ are Markov chain with transition kernel $R$ starting respectively from $x$ and $y$, then we can couple them in a monotone fashion (similar to the coupling of $X$ and $Y$ above), so that
\bq
\PP[\fo n\in \ZZ_+,\, Y_x(n)\leq Y_y(n)]&=&1\eq
(see for instance the book of Lindvall \cite{MR94c:60002}).
\par
Let us prove this monotonicity of $Y$:
\begin{lem}\label{lem5.5}
The Markov chain $Y$ is monotone.
\end{lem}
\proof
Since $Y$ is a birth and death chain, to get it is monotone, it is sufficient to check 
that 
\bq
\fo x\in\lin 0, N-5\rin,\qquad R(x, \lin x-1, x\rin)&\geq & R(x+1,x)\eq
(again see e.g. Lindvall \cite{MR94c:60002}).
\par
The previous bound amounts to 
\bq
\fo x\in\lin 0, N-5\rin,\qquad 1-\f{N-x-1}{N(N+1)}&\geq & \f{(x+1)(N-x-1)}{N(N+1)}\eq
or
\bqn{mon}
\fo x\in\lin 0, N-5\rin,\qquad N(N+1)&\geq &(x+2)(N-x-1)\eqn
\par
The maximum of the r.h.s.\ as $x$ runs in $\RR$ is attained at the point
$x=(N-3)/2$ and replacing in the above r.h.s., the desired inequality is true if we have $N\geq (N+1)/4$, which is satisfied as soon as $N\geq 1/3$.\wwtbp\par
Let us come back to the quantities $\PP[\wi Z(n)>0]$ and $\PP[\wit Z(n)>0]$, we have:
\begin{lem}
For any $n\in\ZZ_+$, we have
\bq
\PP [\wit Z(n)>0]&\leq & \f{2^{N+1}n}{N!}\\
\PP [\wi Z(n)>0]&\leq & \f{ 2^{N+1}n}{N!}\eq
\end{lem}
\proof
We have
\bqn{witZ}
\nonumber\PP[\wit Z(n)>0]&\leq &\EE[\wit Z(n)]\\
&=&\sum_{k=0}^{n-1}\EE[\un_{X(k)\geq Y(k), \,X(k+1)< Y(k+1)}]\eqn
\par
Fix some $k\in\lin 0, n-1\rin$. If $X(k)\geq Y(k)$ and  $X(k+1)< Y(k+1)$ hold, then either $X(k)=Y(k)$ or $X(k)=Y(k)+1$.
Let us consider the latter case, we have:
\bqn{Am}
\EE[\un_{X(k)= Y(k)+1, \,X(k+1)< Y(k+1)}]&=&\sum_{x=0}^{N-4}\PP[X(k)=x+1, Y(k)=x]S((x+1,x), A_-)\eqn
where 
\bq
A_-&\df&\{(x',y')\in\lin 0, N-4\rin\st x'<y'\}\eq
\par
But for the transition from $(x+1,x)$ to $B$ to happen, the underlying uniform random variable on $[0,1]$ must have taken advantage of the discrepancy between $2p(x+1)$ and 1,
otherwise the monotonicity of $Y$ leads to a contradiction. We deduce
\bq
S((x+1,x),A_-)&\leq & \f{\vert 1-2p(x+1)\vert}{N(N-1)}\\
&\leq &  \f1{(N-x-1)!}\f1{N(N-1)}\\
&\leq &  \f1{(N-x+1)!}\eq
and it follows, as in proof of Lemma \ref{lem5.4} that
\bq
\sum_{x=0}^{N-4}\PP[X(k)=x+1, Y(k)=x]S((x+1,x), A_-)&\leq &
\sum_{x=0}^{N-4}\f{1}{Z_N x!}\f1{(N-x+1)!}\\
&\leq & \f{ 2^{N+1}}{(N+1)!}\\
&\leq & \f{ 2^{N}}{N!}\
\eq
\par
The treatment of the cases $X(k)=Y(k)$ is similar to the proof of Lemma \ref{lem5.4}, leading to
\bqn{Am2}
\sum_{x=0}^{N-4}\PP[X(k)=x, Y(k)=x]S((x,x), A_-)&\leq &\f{ 2^{N}}{N!}
 \eqn
\par
It follows that for any $k\in\lin 0,n-1\rin$
\bq
\EE[\un_{X(k)\geq Y(k), \,X(k+1)< Y(k+1)}]&\leq & 
 \f{2^{N+1}}{N!}
\eq
and \eqref{witZ} leads to the first desired bound.\par\me
The second desired bound is obtained in a similar way, the main difference being that we have to replace, for $k\in\lin 0, n-1\rin$, \eqref{Am} by
\bq
\EE[\un_{X(k)= Y(k)-1, \,X(k+1)< Y(k+1)}]&=&\sum_{x=0}^{N-4}\PP[X(k)=x+1, Y(k)=x]S((x+1,x), A_+)\eq
where
\bq
A_+&\df&\{(x',y')\in\lin 0, N-4\rin\st x'>y'\}\eq\par
Then we rather use
\bq
S((x-1,x),A_+)&\leq & \f{\vert 1-2p(x-1)\vert}{N(N-1)}\\
&\leq &  \f1{(N-x-1)!}\f1{N(N-1)}\\
&\leq &  \f1{(N-x+1)!}\eq
leading to 
\bq
\EE[\un_{X(k)= Y(k)-1, \,X(k+1)< Y(k+1)}]&\leq &\f{ 2^{N+1}}{(N+1)!}\\
&\leq &\f{ 2^{N}}{N!}
\eq
\par
As in \eqref{Am2}, we also have
\bq
\EE[\un_{X(k)= Y(k), \,X(k+1)< Y(k+1)}]&\leq &\f{ 2^{N}}{N!}
\eq
enabling us to conclude to the second desired bound.\wwtbp
\par
We are left with the evaluation of the tails of $\tau_0^X$ and $\tau_0^Y$ in \eqref{tails}.\par
We start with the last one:
\begin{lem}\label{lem5.7}There exists a constant $c>0$ such that 
for any $N$ large enough and any $n\in\ZZ_+$, we have
\bq
\PP[\tau_0^Y>n]&\leq & e^{1-cn/N^3}\eq
whatever the initial law of $Y(0)$.
\end{lem}
\proof
For any time $k\in\ZZ_+$ such that $Y(k)=y\neq 0$, we compute
\bq
\EE\lt[\exp\lt.\lt(\frac{Y_{k+1}-Y_k}{N}\rt)\rt\vert Y_k=y\rt]&=&e^{-1/N}\f{y(N-y)}{N(N-1)}+e^{1/N}\f{N-y-1}{N(N-1)}+1-\f{y(N-y)+N-y-1}{N(N-1)}\\
&=&1+(e^{-1/N}-1)\f{y(N-y)}{N(N-1)}+(e^{1/N}-1)\f{N-y-1}{N(N-1)}
\eq
\par
Denoting $F_N(y)$ the r.h.s., it is a second order polynomial whose minimal value is attained at
\bq
\underline{y}&\df& \f{e^{1/N}-1+N}{2(1-e^{-1/N})}
\eq
belonging to $\lin 1, N-4\rin$ for $N$ large enough.
It follows that the maximal value of $F_N(y)$ for $y\in\lin 1, N-4\rin$ is attained either at $y=1$ or $y=N-4$.\par
We compute that 
\bq
F_N(1)&=&1+\f{(e^{-1/N}-1)(N-1)+(e^{1/N}-1)(N-2)}{N(N-1)}\\
F_N(N-4)&=&1+\f{4(e^{-1/N}-1)(N-4)+3(e^{1/N}-1)}{N(N-1)}\eq
and we deduce there exists a constant $c>0$ such that for $N$ large enough,
\bq
 \max(F_N(1),F_N(N-4))&\leq & 1-\f{c}{N^3}\eq
leading to
\bq
\fo k\in\ZZ_+, \qquad Y_k\neq 0&\Rightarrow& \EE\lt[\exp\lt.\lt(\frac{Y_{k+1}}{N}\rt)\rt\vert Y_k\rt]\ \leq \ \lt(1-\f{c}{N^3}\rt)\exp\lt(\frac{Y(k)}{N}\rt)\eq
implying
\bq
\fo k\in\ZZ_+, \qquad \EE\lt[\exp\lt.\lt(\frac{Y_{k+1}}{N}\rt)\rt\vert Y_k\rt]\un_{Y(k)\neq 0}& \leq&\lt(1-\f{c}{N^3}\rt)\exp\lt(\frac{Y(k)}{N}\rt)\eq
i.e.
\bq
\fo k\in\ZZ_+, \qquad \EE\lt[\exp\lt.\lt(\frac{Y_{k+1}}{N}\rt)\un_{Y(k)\neq 0}\rt\vert Y(k)\rt]& \leq&\lt(1-\f{c}{N^3}\rt)\exp\lt(\frac{Y(k)}{N}\rt)\eq
or, using the Markov property,
\bq
\fo k\in\ZZ_+, \qquad \EE\lt[\exp\lt.\lt(\frac{Y_{k+1}}{N}\rt)\un_{Y(k)\neq 0}\rt\vert \cF^Y(k)\rt]& \leq&\lt(1-\f{c}{N^3}\rt)\exp\lt(\frac{Y(k)}{N}\rt)\eq
where $\cF^Y(k)$ is the sigma-field generated by $Y(0), Y(1), ..., Y(k)$.
\par
Iterating this relation, we get for any $k\in\NN$,
\bq
\EE\lt[\lt.\EE\lt[\exp\lt.\lt(\frac{Y_{k+1}}{N}\rt)\un_{Y(k)\neq 0}\rt\vert \cF^Y(k)\rt]\un_{Y(k-1)\neq 0}\rt\vert \cF^Y(k-1)\rt]& \leq&\lt(1-\f{c}{N^3}\rt)^2\exp\lt(\frac{Y(k-1)}{N}\rt)\eq\par
Pushing further the iteration, we end up with
\bq
\lefteqn{\hskip-50mm\EE\lt[\lt.\cdots\EE\lt[\lt.\EE\lt[\exp\lt.\lt(\frac{Y_{k+1}}{N}\rt)\un_{Y(k)\neq 0}\rt\vert \cF^Y(k)\rt]\un_{Y(k-1)\neq 0}\rt\vert \cF^Y(k-1)\rt]\cdots \un_{Y(0)\neq 0}\rt\vert \cF^Y(0)\rt]}\\
&\leq & \lt(1-\f{c}{N^3}\rt)^{k}\exp\lt(\frac{Y(0)}{N}\rt)\eq\par
Taking into account that $Y_{k+1}\geq 0$ and that $Y(0)\leq N$, we get
\bq
\EE\lt[\lt.\cdots\EE\lt[\lt.\EE\lt[\lt.\un_{Y(k)\neq 0}\rt\vert \cF^Y(k)\rt]\un_{Y(k-1)\neq 0}\rt\vert \cF^Y(k-1)\rt]\cdots \un_{Y(0)\neq 0}\rt\vert \cF^Y(0)\rt]
&\leq &e \lt(1-\f{c}{N^3}\rt)^{k}\eq\par
Taking expectation and simplifying conditional expectation iteratively (starting with $\cF^Y(0)$, next $\cF^Y(1)$, etc.), we end up with
\bq
\PP[Y(k)\neq 0, Y(k-1)\neq 0, ..., Y(0)\neq 0]&\leq &e \lt(1-\f{c}{N^3}\rt)^{k}\eq
implying 
\bq
\PP[\tau^Y_0>k]&\leq & e \lt(1-\f{c}{N^3}\rt)^{k}\\
&\leq & e^{1-\f{ck}{N^3}}\eq
which is desired bound, taking $k=n$.\wwtbp
\par
The tail of $\tau_0^X$ is evaluated similarly:
\begin{lem}There exists a constant $\wi c>0$ such that 
for $N$ large enough  and any $n\in\ZZ_+$, we have
\bq
\PP[\tau_0^X>n]&\leq & e^{1-\wi cn/N^3}\eq
whatever the initial law of $X(0)$.
\end{lem}
\proof
According to Proposition \ref{pro4.1}, we have $2p(x)\geq 1/2$ for all $x\in\lin 0, N-4\rin$, fact which suggests to consider Markov chains $\wi Y\df(\wi Y(n))_{n\in\ZZ_+}$ associated
to the transition kernel $\wi R$ given by
\bq
\fo x\neq y\in \lin 0, N-4\rin,\qquad
\wi R(x,y)&\df& \f1{N(N-1)}\lt\{\begin{array}{ll}
x(N-x)&\hbox{, if  $y=x-1$}\\
N-x-1/2&\hbox{, if $y=x+1$}\\
0&\hbox{, otherwise}
\end{array}\rt.\eq
which differs from $R$ only the replacement of $N-x-1$ by $N-x-1/2$.
\par
Consider the corresponding hitting time of 0:
\bq
\tau^{\wi Y}_0&\df&\inf\{n\in\ZZ_+\st \wi Y(n)=0\}\eq
\par Coupling in a monotone way $X$ and $\wi Y$ and starting with $\wi Y(0)=X(0)$, we get 
that
\bq
\fo n\in\ZZ_+, \qquad X(n)&\leq & \wi Y(n)\eq
at least if $\wi Y$ is monotone. This is true and is proven as for Lemma \ref{lem5.5}, where \eqref{mon} has to be replaced by
\bq
\fo x\in\lin 0, N-5\rin,\qquad N(N+1)&\geq &(x+2)(N-x-1/2)\eq
\par
We deduce that
\bq
\fo n\in\ZZ_+,\qquad 
\PP[\tau_0^X>n]&\leq &\PP[\tau_0^{\wi Y}>n]\eq
\par
It is thus sufficient to find 
a constant $\wi c>0$ such that 
for any $N\geq 5$ and any $n\in\ZZ_+$, we have
\bq
\PP[\tau_0^{\wi Y}>n]&\leq & e^{1-\wi cn/N^3}\eq
whatever the initial law of $\wi Y(0)$.
\par
This done as in the proof of Lemma \ref{lem5.7}.\wwtbp
\par
Summarizing the previous computation, we have shown there exist two constants $c,\wi c>0$ such that for any $N$ large enough and $n\geq 0$,
\bq
\lVe \pi-\zeta\rVe_{\mathrm{tv}}&\leq & \f{ 2^{N}n}{N!}+2\f{ 2^{N+1}n}{N!}
+
e^{1- cn/N^3}+e^{1-\wi cn/N^3}\\
&\leq &\f{5\times 2^{N}n}{N!}+2e^{1- \wit cn/N^3}
\eq
with $\wit c\df c\wedge \wi c$.
\par
Taking $n=N\ln(N)/\wit c$, we conclude  \eqref{lim2}.

 \appendix
 
\section{Recovering classical results on $\pi$ through the Markov approach}\label{Api}

Working in the same spirit as in Section \ref{priori}, it is possible to recover the exact formula for the number of fixed points $\pi(x)$ (see \eqref{pi}) from the reversibility of  the Markov chain $P$ (see \eqref{P})   with respect to $\pi$, leading to an alternative proof
for Montmort's formula \eqref{der} (the traditional argument goes through the inclusion-exclusion principle, see e.g.\ \cite{enwiki:1084547992} or Chapter 1 of Arratia, Barbour and Tavaré \cite{zbMATH02042126}).
\par\me
We will use the birth and death chain  $\wit P$ defined in \eqref{witP} above and it's stationary distribution $\wit\pi$ defined in \eqref{witpi}. Using that notation, the reversibility says 
\bq
\fo i\in \lin 0, N-1\rin,\qquad \wit\pi(i)\wit P(i,i+1)&=&\wit \pi(i+1)\wit P(i+1,i)\eq
or
\bq
\fo i\in \lin 0, N-1\rin,\qquad \pi(z_i) P(z_i,z_{i+1})&=& \pi(z_{i+1})P(z_{i+1},z_i)\eq
namely
\bqn{pipi}
\pi(0)P(0,1)\ =\ \pi(1)P(1,0)&\hbox{and}&\fo x\in V\setminus\{N\},\, \pi(x)P(x,x+2)\ =\ \pi(x+2)P(x+2,x)\eqn
i.e.
\bq
\pi(0)(N-2p(0))\ =\ \pi(1)(N-1)&\hbox{and}&\fo x\in V\setminus\{N\},\, 2\pi(x)p(x)\ =\ \pi(x+2)(x+2)(x+1)\eq 
\par
The last condition implies that
\bq
\fo x\in V\setminus\{N\},\qquad p(x)&=& \f{\pi(x+2)}{2\pi(x)}(x+2)(x+1)\eq\par
This formula also holds for $x=N$, since both terms vanish, thus we have shown:\par
\begin{lem}\label{lem4}
We have
\bq\fo x\in V,\qquad p(x)&=& \f{\pi(x+2)}{2\pi(x)}(x+2)(x+1)\eq
\end{lem}
\par
Replacing this expression in the definition of the first associated birth and death kernel  $\wi P$ (defined in \eqref{wiP}), we will deduce the following expression for the reversible probability $\pi$:
\begin{pro}\label{pro1}
We have 
\bq\fo x\in V,\qquad \pi(x)&=&\f1{x!} \sum_{k=0}^{N-x}\f{(-1)^k}{k!}\eq
\end{pro}
\proof
From Lemma \ref{lem4}, we get for any $x\neq y\in V$,
\bq
\wi P(x,y)&=& \f1{N(N-1)}\lt\{\begin{array}{ll}
x(N-x)&\hbox{, if $x\neq N$ and $y=x-1$}\\
N-x-\f{\pi(x+2)}{\pi(x)}(x+2)(x+1)&\hbox{, if $x\neq N-2$ and $y=x+1$}\\
2&\hbox{, if $x= N-2$ and $y=N$}\\
N(N-1)&\hbox{, if $x= N$ and $y=N-2$}\\
0&\hbox{, otherwise}
\end{array}\rt.\eq
\par
Thus for $x\in\lin 0, N-3\rin$, the relation $\pi(x)\wi P(x,x+1)= \pi(x+1)\wi P(x+1,x)$ becomes
\bqn{f1}
\pi(x)(N-x)-\pi(x+2)(x+2)(x+1)&=&\pi(x+1)(x+1)(N-x-1)\eqn
\par
For $x=N-2$, the relation $\pi(N-2)\wi P(N-2,N)= \pi(N)\wi P(N,N-2)$ becomes
\bqn{f2}
\pi(N-2)2&=&\pi(N)N(N-1)\eqn\par
These relations lead us to introduce the function $f$ on $V$ defined by
\bq
\fo x\in \lin 0, N\rin,\qquad f(x)&\df&
\f{\pi(x)}{\cP(x)}\eq
where $\cP$ is the Poisson distribution of parameter 1 (with the convention $f(N-1)=0=\pi(N-1)$).
Indeed, \eqref{f1} and \eqref{f2} reduce to 
\bq
\fo x\in \lin 0, N-3\rin,\qquad f(x)(N-x)-f(x+2)&=&f(x+1)(N-x-1)\\
2f(N-2)&=&f(N)\eq
namely
\bq
\fo x\in \lin 0, N-2\rin,\qquad (f(x)-f(x+1))(N-x)&=&f(x+2)-f(x+1)\eq
\par
This relation leads to the introduction of the function $g$ on $V$ defined by
\bq
\fo x\in\lin 0, N-1\rin,\qquad g(x)&\df&f(x+1)-f(x)\eq
since we get
\bq
\fo x\in\lin 0, N-2\rin,\qquad g(x)&=&-\f{g(x+1)}{N-x}\\
&=&\f{g(x+2)}{(N-x)(N-x-1)}\\
&=&(-1)^{N-x}\f{g(N-1)}{(N-x)!}\\
&=&(-1)^{N-x}\f{f(N)}{(N-x)!}
\eq
\par
Taking into account that $g(N-2)=f(N-1)=0$, we deduce that
\bq
\fo x\in\lin 0, N-2\rin,\qquad
f(x)&=&-g(x)-g(x+1)-\cdots -g(N-2)\\
&=& f(N)\sum_{k=0}^{N-x}\f{(-1)^k}{k!}\eq
a formula also valid for $x=N$, so finally
\bq
\fo x\in V,\qquad\pi(x)&=&\f{1}{A_Nx!} \sum_{k=0}^{N-x}\f{(-1)^k}{k!}\eq
where $A_N=e/f(N)$. This quantity is also the normalization factor, since $\pi$ is a probability,
so we compute
\bq
A_N&=&\sum_{x\in V}\f1{x!}\sum_{k=0}^{N-x}\f{(-1)^k}{k!}\\
&=&\sum_{x=0}^{N-2}\f1{x!}\sum_{k=0}^{N-x}\f{(-1)^k}{k!}+\f1{N!}\\ 
&=&\sum_{x=0}^{N}\f1{x!}\sum_{k=0}^{N-x}\f{(-1)^k}{k!}-\f1{(N-1)!} \sum_{k=0}^{1}\f{(-1)^k}{k!}-\f1{N!}\sum_{k=0}^{0}\f{(-1)^k}{k!}+\f1{N!}\\ 
&=&\sum_{x=0}^{N}\f1{x!}\sum_{l=0}^{N-x}\f{(-1)^{N-x-l}}{(N-x-l)!}-\f1{(N-1)!}(1-1)-\f1{N!}+\f1{N!}\\
&=&\sum_{l=0}^{N}\sum_{x=0}^{N-l}\f{(-1)^{N-x-l}}{x!(N-x-l)!}\\
&=&\sum_{l=0}^{N}\f1{(N-l)!}\lt(1-1\rt)^{N-l}\\
&=&1
\eq
\wwtbp
\par
From this formula, we recover an upper bound on the total variation distance between $\pi$ and $\cP$
almost as good as that of  \eqref{theo1}, but which is not going through a coupling.
Indeed, we compute:
\bq
\lVe \pi-\cP\rVe&=&\sum_{n\in\ZZ_+} (\pi(n)-\cP(n))_+\\
&=&\sum_{n\in\lin 0, N\rin} \lt(\f{D_{N-n}}{(N-n)!}-e^{-1}\rt)_+ \f1{n!}\\
&=&\sum_{n\in\lin 0, N\rin} \lt(\f{D_{n}}{n!}-e^{-1}\rt)_+ \f1{(N-n)!}\\
&=&\sum_{n\in\lin 0, N\rin} \lt( \sum_{k=0}^{n}\f{(-1)^k}{k!}-e^{-1}\rt)_+ \f1{(N-n)!}\\
&=&\sum_{n\in\lin 0, N\rin} \lt( \sum_{k\geq n+1}\f{(-1)^k}{k!}\rt)_+ \f1{(N-n)!}\\
&\leq & \sum_{n\in\lin 0, N\rin,\, n\,\mathrm{odd}}\f{1}{(n+1)!}\f1{(N-n)!}
\eq
(where we used the alternance of the terms of the series $\sum_{k\geq 0}\f{(-1)^k}{k!}$).\par
The last term is also equal to
\bq
 \sum_{n\in\lin 0, N\rin,\, n\,\mathrm{even}}\f{1}{n!}\f1{(N+1-n)!}&\leq & \sum_{n\in\lin 0, N\rin}\f{1}{n!}\f1{(N+1-n)!}\\
 &=&\f1{(N+1)!} \sum_{n\in\lin 0, N\rin}\binom{N+1}{n}\\
 &=&\f{2^{N+1}}{(N+1)!}\eq
\par\sm
For completeness, let us also recall a simple proof of the well-known formula \eqref{piD}:
\begin{lem}\label{lemA.3}
For any $x\in V$, we have 
\bq\pi(x)&=&\f{D_{N-x}}{(N-x)!}\f1{x!}\eq\par
\end{lem}
\proof
Fix $x\in V$ and denote $I_x$ the set of subsets of $\lin N\rin$ whose cardinal is $x$.
By symmetry we have, denoting by $\sigma$ a generic permutation,
\bq
\pi[\eta_1=x]&=&\sum_{I\in I_x}\pi[\fo i\in I,\, \sigma(i)=i,\,\fo j\in\lin N\rin\setminus I,\, \sigma(j)\neq j]\\
&=&\binom{N}{x}\pi[\fo j\in\lin N-x\rin,\, \sigma(j)\neq j,\, \fo i\in \lin N-x+1,N\rin,\, \sigma(i)=i]\\
&=&\binom{N}{x}\f{D_{N-x}}{N!}\\
&=&\f{D_{N-x}}{(N-x)!}\f1{x!}\eq\wwtbp
\par
Montmort's formula \eqref{der} is now a consequence of the above lemma and of  Proposition \ref{pro1}:
\begin{cor}
We have for any $n\in\NN$,
\bq
D_n&=&n!\sum_{k=0}^n\f{(-1)^k}{k!}\eq
\end{cor}
\par
\begin{rem}\par a) It seems from \eqref{pipi} that we have an extra relation for $p(0)$: $\pi(0)(N-2p(0))\ =\ \pi(1)(N-1)$, which amounts to
\bq
p(0)&=&\f{N}2\lt(1-\f{D_{N-1}}{D_N}(N-1)\rt)\eq
\par
Comparing with \eqref{qD}, which gives for $x=0$,
\bq
p(0)&=&\f12\f{D_{N-2}}{(N-2)!}\f{N!}{D_{N}}\eq
we deduce 
\bq D_N&=&(N-1)(D_{N-1}+D_{N-2})\eq
\par
This is the well-known iteration formula for the derangement numbers, see e.g.\ \cite{enwiki:1084547992}.
\par\sm
b) Note that $\pi$ is not close to $\cP$ is the separation discrepancy 
\bq
\fs(\pi,\cP)&=&\sup\lt\{1-\f{\pi(x)}{\cP(x)}\st x\in\ZZ_+\rt\}\eq
since the r.h.s.\ is trivially 1. But with the notations of Section \ref{coupling}, we even  have
\bq
\liminf_{N\ri\iy} \fs(\check \pi,\zeta)&\geq &
\liminf_{N\ri\iy}1-\f{\check\pi(N-4)}{\zeta(N-4)}\\
&=&\lim_{N\ri\iy}1-\f{\pi(N-4)}{\cP(N-4)}\\
&=&1-e\f{D_4}{4!}\\
&>&0\eq
\par
This fact a priori excludes a proof via strong stationary times (see Diaconis and Fill \cite{MR1071805}) in Section~\ref{coupling}.
\end{rem}
\par

\section{Complements on the conditional expectation $p$}\label{complements}

Some observations about $p$ are gathered here.
\par\me
Note that Lemma \ref{lemA.3} also leads to an expression of the quantities $p(x)$ in terms of the number of derangements, from Lemma \ref{lem4}:
\bqn{qD}
\fo x\in V,\qquad p(x)&=&\f12\f{D_{N-x-2}}{(N-x-2)!}\f{(N-x)!}{D_{N-x}}\\
\nonumber&=&\f12\f{\sum_{k=0}^{N-x-2}\f{(-1)^k}{k!}}{\sum_{l=0}^{N-x}\f{(-1)^l}{l!}}
\eqn
\par
This formula leads to an estimate of our quantities of interest, the $\vert 2p(x)-1\vert $, for $x\in\lin 0, N-2\rin$, of the same order as that of Proposition \ref{pro4.1}:
\begin{lem}\label{lemB1}
We have
\bq
\fo x\in\lin 0, N-2\rin,\qquad \vert 2p(x)-1\vert &\leq & 3\f{N-x-1}{(N-x)!}\\
&\leq & \f3{(N-x-1)!}
\eq and in particular we get, for $N\geq 4$,
\bq
\fo x\in\lin 0, N-4\rin,\qquad \f14\ \leq \ p(x)\ \leq \ \f34\eq
\end{lem}
\proof
From \eqref{qD} we deduce:
\bq
\fo x\in V,\qquad 2p(x)&=&\f{\sum_{k=0}^{N-x-2}\f{(-1)^k}{k!}}{\sum_{l=0}^{N-x}\f{(-1)^l}{l!}}\\
&=&1-\f{\sum_{N-x-1}^{N-x}\f{(-1)^k}{k!}}{\sum_{l=0}^{N-x}\f{(-1)^l}{l!}}
\eq
implying
\bq
\fo x\in V,\qquad \vert 2p(x)-1\vert &=&
\f{\lve\f{1}{(N-x-1)!}-\f1{(N-x)!}\rve}{\sum_{l=0}^{N-x}\f{(-1)^l}{l!}}\\
&=&\f{1}{(N-x-1)!}\f{1-\f1{N-x}}{\sum_{l=0}^{N-x}\f{(-1)^l}{l!}}\\
&=&\f{N-x-1}{(N-x)!}\f1{\sum_{l=0}^{N-x}\f{(-1)^l}{l!}}
\eq
\par
Note that the series $\sum_{l=0}^{n}\f{(-1)^l}{l!}$ provide alternating approximations of $e^{-1}$, it follows that
\bq
\fo x\in\lin 0, N-2\rin,\qquad \sum_{l=0}^{3}\f{(-1)^l}{l!}\ \leq \ \sum_{l=0}^{N-x}\f{(-1)^l}{l!}\ \leq \ \sum_{l=0}^{2}\f{(-1)^l}{l!}\eq
namely
\bq
\fo x\in\lin 0, N-2\rin,\qquad \f1{2!}-\f1{3!}\ \leq \ \sum_{l=0}^{N-x}\f{(-1)^l}{l!}\ \leq \ \f1{2!}\eq
i.e.
\bq
\fo x\in\lin 0, N-2\rin,\qquad \f1{3} \leq \ \sum_{l=0}^{N-x}\f{(-1)^l}{l!}\ \leq \ \f1{2}\eq
whose lower bound leads to the first desired estimate.\par
For the second estimate, note that
\bq
\fo x\in\lin 0,N-4\rin,\qquad
\f3{(N-x-1)!}&\leq & \f3{(N-(N-4)-1)!}\\
&\leq &\f3{3!}\ =\ \f12\eq
\qquad
\wwtbp
\par
Lemma \ref{lemB1} can be used similarly to Proposition \ref{pro4.1} in Section \ref{coupling}, leading to the same conclusion.\par\me
Coming back to the formulation \eqref{p} of $p$ as a conditional expectation of $\eta_2$, it is natural to wonder if it could not be deduced from symmetry arguments.
Remark it is true for the whole expectation: $\EE_\nu[\eta_2]=1/2$ (see the proof of Lemma \ref{lem19} below with $k=0$), in the same way one immediately gets $\EE_\nu[\eta_1]=1$. So to finish this appendix, let us 
show that symmetry arguments lead  to a natural linear equation satisfied by $p$, even if we did not find how to use it to deduce the a priori bounds similar to those of Proposition \ref{pro4.1} or Lemma \ref{lemB1}.\par\sm
For $k\in \lin 0, N\rin$, denote 
\bq
\cA_k&=&\{(i_1, i_2, ..., i_k)\in \lin N\rin^k\st m\neq n\in \lin k\rin\Rightarrow i_m\neq i_n\}\eq
\par
In particular, we have
\bqn{vAk}
\vert \cA_k\vert&=& N(N-1)\cdots (N-k+1)\eqn
(by convention, $\cA_0=\{\emptyset\}$ and $\vert \cA_0\vert=1$).\par
For $k\in\lin0, N\rin$, we define the mapping $F_k$ on the symmetric group $\cS_N$ via
\bq
\fo \sigma\in\cS_N,\qquad F_k(\sigma)&\df& \sum_{(i_1, ..., i_k)\in \cA_k} \prod_{j\in\lin k\rin} \un_{\{\sigma(i_j)=i_j\}}\eq
\par
Let us check these mappings are functions of $\eta_1$ (the number of fixed points):
\begin{lem}
For any $k\in\lin 0, N\rin$, we have
\bqn{Feta}
F_k&=&
\eta_1(\eta_1-1)\cdots (\eta_1-k+1)\eqn
\end{lem}
\proof
Indeed, for any given $\sigma\in\cS_N$, denote $\cF(\sigma)$ the set of fixed points of $\sigma$. We have
\bq
F_k(\sigma)&=&\vert \cA_k\cap \cF(\sigma)^k\vert\\
&=&\eta_1(\sigma)(\eta_1(\sigma)-1)\cdots (\eta_1(\sigma)-k+1)\eq
\wwtbp
\par
\begin{rem}\label{poly}
Since $F_k$ is a polynomial of order $k$ in $\eta_1$, any function of $\eta_1$ can be expressed as a linear combination of the $F_k$ for $k\in\lin 0, N\rin$, and even only for $k\in\lin 0, N-1\rin$ or alternatively $k\in V$, since $\eta_1$ is taking $N$ values, those of $V\df\lin 0, N\rin\setminus\{N-1\}$. 
\end{rem}
\par
It follows that if we want to prove that 
\bq
\EE_\nu[ \eta_2\vert \eta_1]&=&f(\eta_1)\eq
for a given function $f\st V\ri \RR_+$, it is sufficient to check that
\bq
\fo k\in\lin 0, N\rin,\qquad \EE_\nu[ \eta_2 F_k]&=&\EE_\nu[ f(\eta_1) F_k]\eq
\par
We are thus led to compute the l.h.s.\par
\begin{lem}\label{lem19}
For any $k\in\lin 0, N\rin$, we have
\bq
\EE_\nu[ \eta_2 F_k]&=&\lt\{\begin{array}{ll}
1/2&\hbox{, if $k\in\lin 0, N-2\rin$}\\[2mm]
0&\hbox{, if $k\in\{N-1,N\}$}\end{array}\rt.
\eq
\end{lem}
\proof
Note that
\bq
\fo \sigma\in\cS_N,\qquad
\eta_2(\sigma)&=&\f12\sum_{m\in \lin N\rin}\un_{\{\sigma(m)\neq m,\, \sigma^2(m)=m\}}\\
&=&\f12\sum_{m\neq n\in \lin N\rin}\un_{\{\sigma(m)=n,\, \sigma(n)=m\}}
\eq
so that
\bq
2\EE_\nu[ \eta_2 F_k]&=&\sum_{m\neq n\in \lin N\rin}\sum_{(i_1, ..., i_k)\in \cA_k}\EE_\nu\lt[ \un_{\{\sigma(m)=n,\, \sigma(n)=m\}}\prod_{j\in\lin k\rin} \un_{\{\sigma(i_j)=i_j\}}\rt]
\eq
\par
Note that the above expectation vanishes if $m\in\{i_1, ..., i_k\}$ or $n\in\{i_1, ..., i_k\}$, so writing $i_{k+1}=m$ and $i_{k+2}=n$, we end up with
\bq
2\EE_\nu[ \eta_2 F_k]&=&\sum_{(i_1, ..., i_k,i_{k+1},i_{k+2})\in \cA_{k+2}}\EE_\nu\lt[ \un_{\{\sigma(i_{k+1})=i_{k+2},\, \sigma(i_{k+2})=i_{k+1}\}}\prod_{j\in\lin k\rin} \un_{\{\sigma(i_j)=i_j\}}\rt]\\
&=&\sum_{(i_1, ..., i_k,i_{k+1},i_{k+2})\in \cA_{k+2}}\!\!\!\!\pi[\sigma(i_1)=i_1, \sigma(i_2)=i_2, ..., \sigma(i_k)=i_k,\sigma(i_{k+1})=i_{k+2},\, \sigma(i_{k+2})=i_{k+1}]
\eq
\par
For any $(i_1, ..., i_k,i_{k+1},i_{k+2})\in \cA_{k+2}$, the above probability can be computed by first choosing $\sigma(i_1)=i_1$, whose probability is $1/N$, next choosing 
$\sigma(i_2)=i_2$, whose subsequent probability is $1/(N-1)$, etc, up to choosing $\sigma(i_{k+2})=i_{k+1}$, whose probability is $1/(N-k-1)$.
We deduce 
\bq
\pi[\sigma(i_1)=i_1, \sigma(i_2)=i_2, ..., \sigma(i_k)=i_k,\sigma(i_{k+1})=i_{k+2},\, \sigma(i_{k+2})=i_{k+1}]&=&\f{1}{N(N-1)\cdots (N-k-1)}\eq
and by consequence
\bq
2\EE_\nu[ \eta_2 F_k]&=&\f{\vert \cA_{k+2}\vert}{N(N-1)\cdots (N-k-1)}\\
&=&1\eq
due to \eqref{vAk}, at least when $k+2\leq N$. Obviously, when $k\in\lin 0, N-1\rin$ satisfies $k+2>N$, namely when $k\in\{N-1,N\}$, we end up with
$2\EE_\nu[ \eta_2 F_k]=0$.\par\wwtbp
\par
\begin{lem}\label{lemB.3}
For any $k\in\lin 0, N\rin$, we have
\bq
\EE_\nu[  F_k]&=&1
\eq
\end{lem}
\proof
Indeed, as in the above proof,
\bq
\EE_\nu[F_k]&=&\sum_{(i_1, ..., i_k)\in \cA_k}\pi[\sigma(i_1)=i_1, \sigma(i_2)=i_2, ..., \sigma(i_k)=i_k]\\
&=&\sum_{(i_1, ..., i_k)\in \cA_k}\f{1}{N(N-1)\cdots (N-k+1)}\\
&=& \f{\vert \cA_{k}\vert}{N(N-1)\cdots (N-k+1)}\\
&=&1
\eq
\wwtbp
\par
It follows that if $f\st V\ri \RR$ is a function satisfying
\bqn{cond}
\fo k\in V,\qquad \EE_\nu[f(\eta_1)F_k]&=&1\eqn
then we can conclude that $f=\un$, the function only taking the value 1.
\par
Consider $\st V\ri \RR$ given by the conditional expectation
\bq
f(\eta_1)&=&\EE[2\eta_2\vert \eta_1]
\eq
\par
According to Lemma \ref{lem19}, $f$ almost satisfies \eqref{cond}, the only discrepancy being the case $k=N$.
Of course it can not satisfy \eqref{cond}, otherwise we would get from Section \ref{coupling} that the law of $\eta_1$ is the conditioning of the Poisson distribution of parameter 1 to $V$ and this is not true (e.g.\ due to \eqref{pi}).\par
\me
Nevertheless, Lemma \ref{lem19} leads to a linear equation for $f$. Denote $a\df(a_k)_{k\in V}$ the vector of the coefficients in the writing
\bq
f(\eta_1)&\fd&\sum_{k\in V} a_k F_k\eq
we have
\bq
G a&=&\lt(\begin{array}{c} 1\\1\\ \vdots\\ 1\\ 0\end{array}\rt)\eq
i.e.
\bq
a&=&G^{-1} \lt(\begin{array}{c} 1\\1\\ \vdots\\ 1\\ 0\end{array}\rt)\eq
where $G\df(G_{k,l})_{k,l\in V}$ is the Gram matrix given by
\bqn{G}
\fo k,l\in V,\qquad G_{k,l}&\df& \EE_\nu[F_kF_l]\eqn
\par
In accordance with Remark \ref{poly}, the family $(F_k)_{k\in V}$ is linearly independent in $\LL^2(\pi_1)$, due to the fact that $\pi_1(x)>0$ for all $x\in V$, which implies that $\dim(\LL^2(\pi_1))=N$. As a consequence, $G$ is invertible.\par
As seen in Section \ref{coupling}, more interesting for us is the function $g\df f-\un$ defined on $V$.
Since Lemma \ref{lemB.3} shows that
\bq
\un&=&\sum_{k\in V}  b_k F_k\eq
with
\bq
\lt(\begin{array}{c} b_0\\b_1\\ \vdots\\ b_{N-2}\\ b_{N-1}\end{array}\rt)&=&G^{-1} \lt(\begin{array}{c} 1\\1\\ \vdots\\ 1\\ 1\end{array}\rt)\eq
we deduce that $g(\eta_1)=\sum_{k\in V} c_kF_k$, with
\bq
\lt(\begin{array}{c} c_0\\c_1\\ \vdots\\ c_{N-2}\\ c_{N-1}\end{array}\rt)&=&G^{-1} \lt(\begin{array}{c} 0\\0\\ \vdots\\ 0\\ 1\end{array}\rt)\eq
\par
More precisely, the computations of  Section \ref{coupling} show the only a priori informations we need to control our coupling constructions are
 estimates on expressions such as
\bqn{needed}
\sum_{x\in \lin 0, N-2\rin} \vert g(x)\vert \,\f1{e x!}\eqn
\par
Below we compute the entries of $G$  directly via symmetry arguments, without a priori knowledge of the law $\pi$ of $\eta_1$, 
nevertheless, it does not seem very helpful to estimate expressions such as \eqref{needed}. 
\begin{pro}
The matrix $G$ is symmetric and extending Definition \eqref{G} to any $k,l\in \lin 0, N\rin$, we have
\bq
\fo k\leq l\in \lin 0, N\rin, \qquad G_{k,l}&=&k!\sum_{r=0}^{k\wedge (N-l)}\f1{r!}\binom{l}{k-r}
\eq
\end{pro}
\proof
For any $i=(i_1, ..., i_k)\in\cA_k$, denote $\{i\}$ the set $\{i_1, ..., i_k\}\subset \lin N\rin$, as well as
\bq
\wi\cA_k&\df&\{\{i\}\st i\in\cA_k\}\\
&=&\{S\subset \lin N\rin\st \vert S\vert =k\}\eq
\par
We compute, for any $k\leq l\in\lin N\rin$, 
\bq
\EE_\nu[F_kF_l]&=&k!l!\sum_{S\in\wi\cA_k,\, T\in\wi\cA_l} \PP[\fo s\in S, \sigma(s)=s,\,\fo t\in T,\sigma(t)=t]\\
&=&k!l!\sum_{S\in\wi\cA_k,\, T\in\wi\cA_l}\f1{N(N-1)\cdots (N-\vert S\cup T\vert+1)}\\
&=&k!l!\sum_{u\in\lin l, (l+k)\wedge N\rin}{A_{k,l}(u)}\f{(N-u)!}{N!}\eq
where
\bq
\fo u\in\lin l, (l+k)\wedge N\rin, \qquad A_{k,l}(u)&\df& \vert \{(S,T)\in \wi\cA_k\times\wi\cA_l\st \vert S\cup T\vert=u\}\vert\eq\par
\par
Note that for any fixed $T\in\wi \cA_l$ and $r\in \lin 0, k\wedge (N-l)\rin$, we have
\bq
\vert\{S\in\wi \cA_k\st \vert S\setminus T\vert = r\}\vert&=&\binom{N-l}{r}\binom{l}{k-r}\eq
the r.h.s.\ corresponding to the number of choices of $r$ elements in $\lin N\rin\setminus T$ and $k-r$ elements in $T$.\par
It follows that, with the change of variable $u=l+r$,
\bq
A_{k,l}(l+r)&=&\binom{N}{l}\binom{N-l}{r}\binom{l}{k-r}\eq
and by consequence
\bq
G_{k,l}&=&k!l!\sum_{r\in \lin 0, k\wedge (N-l)\rin}\binom{N}{l}\binom{N-l}{r}\binom{l}{k-r}\f{(N-l-r)!}{N!}\\
&=&\f{k!l!}{N!}\binom{N}{l}\sum_{r\in \lin 0, k\wedge (N-l)\rin}\binom{N-l}{r}\binom{l}{k-r}(N-l-r)!\\
&=&\f{k!l!}{(N-l)!}\sum_{r\in \lin 0, k\wedge (N-l)\rin}\f{(N-l)!}{r!(N-l-r)!}\f{l!}{(l-k+r)!(k-r)!}(N-l-r)!\\
&=&k!\sum_{r\in \lin 0, k\wedge (N-l)\rin}\f{1}{r!}\f{l!}{(l-k+r)!(k-r)!}\\
&=&k!\sum_{r=0}^{k\wedge (N-l)}\f1{r!}\binom{l}{k-r}
\eq
\wwtbp

\vskip2cm
\hskip70mm
\vbox{
\copy5
 \vskip5mm
 \copy6
}

\end{document}